\documentclass[letterpaper, 10 pt, conference]{ieeeconf}  

\IEEEoverridecommandlockouts                              
\usepackage{cite}
\usepackage{amsmath,amssymb,amsfonts}
\usepackage{algorithmic}
\usepackage{graphicx}
\usepackage{textcomp}
\usepackage{stfloats} 
\usepackage{xcolor}
   
\newtheorem{theorem}{Theorem}
\newtheorem{corollary}{Corollary}
\newtheorem{remark}{Remark}

\newtheorem{lemma}{Lemma}
\newtheorem{assumption}{Assumption}
\newtheorem{proposition}{Proposition}

\title{\LARGE \bf
Combining Exponential and Polynomial Convergence Rates for Accelerated Extremum Seeking: A High-Order Lie Bracket Approach}

\author{Victoria Grushkovskaya$^{1,3}$ and Sameh A. Eisa$^{2}$
\thanks{$^{1}$Victoria Grushkovskaya is with the Department of Mathematics,
        University of Klagenfurt, 9020 Klagenfurt am W\"orthersee, Austria
        {\tt\small viktoriia.grushkovska@aau.at}}%
\thanks{$^{2}$Sameh A. Eisa is with the Department of Aerospace Engineering and Engineering Mechanics, University of Cincinnati, Ohio, United States {\tt\small eisash@ucmail.uc.edu}}
\thanks{$^{3}$Institute of Applied Mathematics \& Mechanics, National Academy of Sciences of Ukraine}
}

\begin{document}

\maketitle
\thispagestyle{empty}
\pagestyle{empty}

\begin{abstract}
In this paper, we leverage the recently developed general higher-order Lie
bracket approximation theory and framework for control-affine systems to show that we can design extremum seeking (ES) systems
that combine flows based on: (i) first-order Lie bracket approximation, which provide an exponential convergence rate for quadratic objective functions and polynomial convergence rate for non-quadratic objective functions; and (ii) higher-order
 Lie bracket approximation, which provide exponential convergence
rate near the minimizer even if the objective function is non-quadratic (i.e., flatter and possesses vanishing Hessian at the minimizer). Such combination provides an accelerated ES design that benefits from the polynomial decay rate farther away from the minimizer and exponential decay rate near the minimizer even for complicated nature of objective functions. We
provide numerical simulations to demonstrate the results, including multi-variable and non-polynomial flat-bottomed objective functions.  
\end{abstract}

\section{Introduction}
Consider the control-affine system
\begin{equation}
\dot{x}=f_0(x)+\sum_{i=1}^{q}\epsilon^{-\vartheta_i}f_i(x)u_i(k_i t/\epsilon),
\end{equation}
where $x\in\mathbb{R}^n$ denotes the state, $\vartheta_i\in(0,1)$ with $0<\epsilon\ll 1$, $f_0$ is the drift vector field, $f_i$ are control vector fields, $u_i$ are control inputs, $q$ is the number of inputs, and $k_i\in\mathbb{Q}_+$. Systems of this type arise in many fields and application areas, including, for example, robotics, multi-agent systems, and flight dynamics (e.g., \cite{Ar03,Tan10,scheinker2024100,grushkovskaya2018family,ECC_2024,eisa2023,palanikumar2026model,elgohary2025hovering,elgohary2025model,palanikumar2026model2,scheinker2017modelfreebook}). Their structure makes them especially suitable for geometric control techniques based on Lie bracket approximations and relevant model-free optimization and control via extremum seeking \cite{pokhrel2026higher}.

Extremum seeking (ES) is a real-time, model-free optimization and control methodology aimed at steering a dynamical system toward the minimizer of an objective function that can be measured online, even when its analytical expression is unavailable \cite{Ar03,scheinker2024100,Tan10}. Classical ES schemes, many of which rely on averaging methods \cite{Maggia2020higherOrderAvg}, have been widely studied for both analysis and controller design \cite{Kr00,yilmaz2025exponential,guay2015time}. Another line of research (the focus of this paper) exploits control-affine formulations and Lie bracket approximations to analyze ES dynamics. In particular, when $\vartheta_i=0.5$, $k_i=1$ for all $i$, first-order Lie bracket systems have been used to explain the behavior of ES algorithms, whose averaged dynamics resemble a gradient flow of the objective -- see for example \cite{DURR2013,GZE18}. Subsequent studies extended this viewpoint beyond first-order Lie brackets. In particular, second-order Lie bracket approximation made it possible to capture information associated with second derivatives of the cost function, thus connecting ES dynamics with Newton-like behavior \cite{labar2019newton}. More recently, a general framework for constructing higher-order Lie bracket approximations has been developed \cite{pokhrel2026higher}. These approximations were shown to possess higher-order averaging dynamics, which guarantee the closeness of the original trajectories and the corresponding approximate dynamics given small enough $\epsilon$. This development opened the door to ES designs based on higher-order Lie brackets. Numerical evidence has suggested that such designs may converge faster than conventional first-order approaches and even outperform Newton-based schemes \cite{pokhrel2026higher}. Most recently, ES laws were developed that ensure motion of the ES system along higher-order Lie bracket directions \cite{grushkovskaya2025extremum,grushkovskaya2026extremum}, yielding exponential convergence not only for quadratic objective functions but also for non-quadratic, flat-bottomed objective functions with a singular Hessian at the minimum. In particular, the ES design in \cite{grushkovskaya2026extremum} is provided for cost functions $J(x)$ that behave locally near the extremum point $x^*$ as $J(x)-J(x^*)\sim \|x-x^*\|^{m}$ with $ m \geq 2$. Thus, the local behavior need not be dominated by a quadratic term ($m=2$), but may instead be dominated by a higher-degree term ($m>2$). The reader can refer to~\cite{grushkovskaya2026extremum} for numerical examples  that involved fourth-, sixth, and even eight-order degree objectives. It is worth noting that Newton-based and exponential convergence results for ES methods typically rely on non-degenerate quadratic behavior of the objective function near the minimizer, or equivalently on suitable non-zero Hessian conditions, see for example  \cite{espitia2026prescribed,yilmaz2025exponential,yilmaz2026asymptotic}.
%

It is also worth highlighting that ES methods for non-quadratic functions (i.e., $m>2$) maintain a polynomial decay rate if they are gradient-based \cite{GZE18}. 
On the other hand, Newton-based methods may face difficulties when inversion of the Hessian is required, since the Hessian may become singular at the minimizer for such non-quadratic, flat-bottomed objective functions. Moreover, exponential convergence cannot be guaranteed in general in Newton-based methods for non-quadratic objective functions even if only estimation of the Hessian is required. 
%

\textbf{Motivation \& Contribution.} 
While polynomial convergence is generally undesirable near the minimizer $x^\ast$ in optimization methods in general and ES methods in particular, it may be associated with a stronger decrease when the state is farther from $x^\ast$. In contrast, exponential convergence becomes preferable in a neighborhood of the minimizer, where it provides a faster decay \cite{grushkovskaya2026extremum}. This naturally raises the question: can one combine the two mechanisms so as to benefit from the stronger decrease farther from $x^\ast$ while retaining exponential convergence near $x^\ast$? 


In this paper, we leverage the general higher-order Lie bracket approximation framework and its application to ES methods \cite{pokhrel2026higher,grushkovskaya2026extremum} to design ES systems that simultaneously excite Lie brackets of different orders. The resulting Lie bracket dynamics combine the gradient term with higher-order derivative terms, benefiting from the stronger polynomial-type decrease away from the minimizer and exponential decay near the minimizer whenever a linear-like higher-order term is present. Hence, we obtain an accelerated ES design based on first- and higher-order Lie brackets. Moreover, the proposed design supports and sheds the light on future possibility where accelerated and exponentially convergent ES becomes possible even if the local degree of the objective function is unknown (i.e., $m$ is totally unknown) since we demonstrate that Lie brackets of different orders can be excited simultaneously. We provide theoretical convergence results and numerical simulations, including multi-variable and non-polynomial flat-bottomed objective functions, to illustrate the proposed approach.

\section{Problem Statement and Main Idea}

We consider a class of extremum seeking (ES) problems for unconstrained minimization of a single-variable cost function $J:D\to\mathbb{R}$, $D\subset\mathbb{R}$:
$$
\min_{x\in D}J(x).
$$
The function $J$ is assumed to be unknown in terms of an explicit analytic expression, but it can be evaluated (measured) at any point. We assume that $J$ possesses an isolated minimum at some $x^\ast\in D$. 



The main idea is to simultaneously excite Lie brackets \cite{grushkovskaya2026extremum} of several odd orders and construct an extremum seeking system whose associated Lie bracket system has the form
\begin{equation}
\label{Lie_sys}
\dot{\bar{x}}=-\sum_{k=0}^{q}\gamma_{2k+1}J^{(2k+1)}(\bar{x}),\, \gamma_{2k+1}>0,
\end{equation}
where $q\in\mathbb{N}\cup\{0\}$. This form combines the gradient direction with directions generated by higher-order derivatives of the cost function.

To illustrate the resulting behavior, suppose that
$
J(x)-J(x^\ast)\sim |x-x^\ast|^m
$
near $x^\ast$. Then $J'(x)$ behaves as $(x-x^\ast)^{m-1}$, whereas $J^{(2k+1)}(x)$ behaves as $(x-x^\ast)^{m-2k-1}$. Hence, the gradient term may provide a stronger decrease farther from $x^\ast$, while the higher-order derivatives become increasingly important close to the minimum. In particular, if $J^{(m-1)}$ is included in~\eqref{Lie_sys}, it generates a linear-like term near $x^\ast$ and therefore exponential decay. Thus, the combination in~\eqref{Lie_sys} preserves the stronger gradient action farther from the minimum while improving the local decay through higher-order derivative terms.

\begin{figure}[t]
    \centering
    \includegraphics[width=1\linewidth]{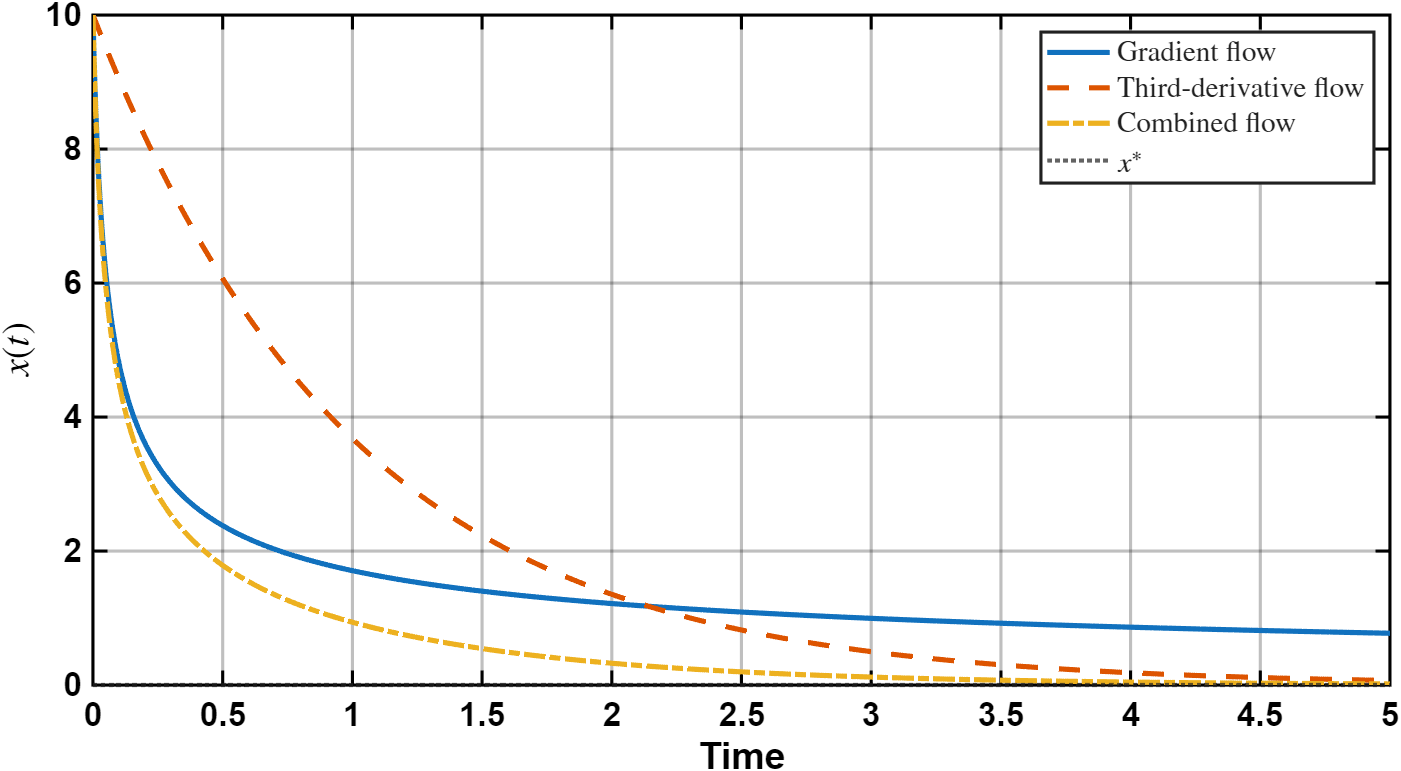}
    \caption{Comparison of the gradient, third-derivative, and combined flows for $J(x)=\frac{1}{4!}(x-x^\ast)^4$, $x^\ast=0$.}
    \label{fig_simple_example}
\end{figure}

For example, let $J(x)=\frac{1}{4!}(x-x^\ast)^4$. Using the first- and third-order derivatives gives the Lie bracket system
$$
\dot{\bar{x}}=-\frac{\gamma_1}{3!}(\bar{x}-x^\ast)^3-\gamma_3(\bar{x}-x^\ast).
$$
The cubic term is dominant for $|\bar{x}-x^\ast|>\sqrt{6\gamma_3/\gamma_1}$, while the linear term is dominant close to $x^\ast$ and provides exponential decay. Thus, the combination takes advantage of the stronger gradient action farther from the minimum and the linear-like behavior generated by the higher-order derivative near $x^\ast$, see Fig.~\ref{fig_simple_example}.

Exciting Lie brackets of several orders may also help with wrong assumptions on the local degree of $J$.  Indeed, choosing a single high-order Lie bracket \cite{grushkovskaya2026extremum} based on an incorrect estimate of the local degree may lead to an ineffective ES law. For example, if $J(x)=\frac{1}{2}(x-x^\ast)^2$ and only a third-order derivative direction is generated, then $J'''(x)\equiv0$ and the corresponding Lie bracket system has zero dynamics. Thus, a design based only on a prescribed higher-order derivative may lose the desired convergence properties when the actual local degree of the cost function is lower than expected. By exciting several odd-order Lie brackets simultaneously, the lower-order derivative terms remain present and provide the required descent direction.

\section{Main Results}

\subsection{Stability Properties of the  Lie Bracket System}

We first analyze the stability and decay properties of the  Lie bracket system~\eqref{Lie_sys}. For this purpose, we introduce the following assumption on the cost function.

\begin{assumption}
\textit{The function $J\in C^{2q+1}(D;\mathbb{R})$ satisfies the following properties:\\
A1.1) there exists an $x^\ast\in D$ such that $J(x^\ast)=J^\ast$, $J(x)>J^\ast$ for all $x\in D\setminus\{x^\ast\}$;\\
A1.2) for every $k\in\{0,\ldots,q\}$, there exist constants $\beta_{k1},\beta_{k2}\geq 0$ and $p_k\geq 1$ such that, for all $x\in D$,
$$
\begin{aligned}
&J^{(2k+1)}(x)(x-x^\ast)\geq \beta_{k1}|x-x^\ast|^{p_k+1},\\
&|J^{(2k+1)}(x)|\leq \beta_{k2}|x-x^\ast|^{p_k},
\end{aligned}
$$
and there exists $\bar k\in\{0,\dots,q\}$ such that $\beta_{\bar k1}>0$.}
\end{assumption}

Assumption A1.2) requires that each excited odd-order derivative has a nonnegative contribution in the direction of the minimum and satisfies a corresponding growth bound. The exponents $p_k$ characterize the asymptotic decay rate associated with the non-vanishing terms: smaller values of $p_k$ correspond to faster asymptotic decay, while $p_k=1$ gives a linear-like term near $x^\ast$.
For sufficiently smooth cost functions whose local Taylor expansion around $x^\ast$ has a finite-order leading term, the nonzero relevant derivatives typically satisfy
$$
J^{(2k+1)}(x)=c_k(x-x^\ast)^{p_k}+o\left(|x-x^\ast|^{p_k}\right),
$$
with $c_k>0$ and odd $p_k$. This local behavior is not restricted to polynomial objective functions. More generally, Assumption~A1.2) only imposes the corresponding sign and growth conditions and does not require a particular Taylor representation.

Similarly to the results in~\cite{GZE18,grushkovskaya2026extremum}, we can characterize the convergence properties of the solutions of~\eqref{Lie_sys} as follows. Define
$$
p^\ast:=\min\limits_{k=0,\ldots,q}\{p_k:\beta_{k1}>0\},
\,
\lambda^\ast:=\sum\limits_{\substack{k=0,\ldots,q\\ p_k=p^\ast,\ \beta_{k1}>0}}\gamma_{2k+1}\beta_{k1}.
$$
\begin{lemma}\label{lemma_lbs_decay}
\textit{If Assumption~1 holds, then $x^\ast$ is an asymptotically stable equilibrium of~\eqref{Lie_sys}. Namely, for any $\delta>0$ such that $B_\delta(x^\ast)\subset D$ and any $\bar{x}(0)\in B_\delta(x^\ast)$, the following estimates hold: 
\begin{itemize}
    \item If $p^\ast=1$, then $x^\ast$ is exponentially stable, i.e., $$
|\bar{x}(t)-x^\ast|\leq|\bar{x}(0)-x^\ast|e^{-\lambda^\ast t}.
$$
\item  If $p^\ast>1$, then $x^\ast$ is polynomially  stable, i.e.,
$$
|\bar{x}(t)-x^\ast|\leq\left(|\bar{x}(0)-x^\ast|^{1-p^\ast}+(p^\ast-1)\lambda^\ast t\right)^{-\frac{1}{p^\ast-1}}.
$$
\end{itemize}}
\end{lemma}

\begin{proof}
Consider $V(\bar{x})=\frac{1}{2}|\bar{x}-x^\ast|^2$. From Assumption~1,
$$
\dot V\leq-\sum_{k=0}^{q}\gamma_{2k+1}\beta_{k1}|\bar{x}-x^\ast|^{p_k+1}\leq-\lambda^\ast|\bar{x}-x^\ast|^{p^\ast+1}.
$$
Thus, $V$ is decreasing, and every solution starting in $B_\delta(x^\ast)$ remains in $B_\delta(x^\ast)$. Moreover, the above estimate implies asymptotic stability of $x^\ast$.

If $p^\ast=1$, then $\dot V\leq-2\lambda^\ast V$, which yields
$
V(t)\leq V(0)e^{-2\lambda^\ast t},
$
and hence
$
|\bar{x}(t)-x^\ast|\leq|\bar{x}(0)-x^\ast|e^{-\lambda^\ast t}.
$

For $p^\ast>1$, we have
$
\dot V\leq-2^{\frac{p^\ast+1}{2}}\lambda^\ast V^{\frac{p^\ast+1}{2}}.
$
Integration of this inequality and substitution of $V=\frac{1}{2}|\bar{x}-x^\ast|^2$ give
$$
|\bar{x}(t)-x^\ast|\leq\left(|\bar{x}(0)-x^\ast|^{1-p^\ast}+(p^\ast-1)\lambda^\ast t\right)^{-\frac{1}{p^\ast-1}}.
$$
\end{proof}
Lemma~1 characterizes the asymptotic decay rate of the solutions of~\eqref{Lie_sys}, which is determined by the smallest exponent $p^\ast$. However, the derivative terms corresponding to larger exponents may have a stronger effect when the state is farther from $x^\ast$. To characterize this property, define
$$
\hat p:=\max\limits_{k=0,\ldots,q}\{p_k:\beta_{k1}>0\},\,
\hat\lambda:=\sum\limits_{\substack{k=0,\ldots,q\\ p_k=\hat p,\ \beta_{k1}>0}}\gamma_{2k+1}\beta_{k1}.
$$
If $p^\ast=1$ and $\hat p>1$, the Lie bracket system contains both a linear-like term determining the asymptotic exponential decay and higher-power terms which may provide a faster decrease farther from the minimum. The following result gives an estimate which takes both effects into account.
\begin{theorem}\label{thm_lbs_decay}
\textit{Let Assumption~1 hold and suppose that $p^\ast=1$ and $\hat p>1$. Then, for any $\delta>0$ such that $B_\delta(x^\ast)\subset D$ and any $\bar{x}(0)\in B_\delta(x^\ast)$,
$$
\begin{aligned}
|\bar{x}(t)-x^\ast|\leq|\bar{x}(0)&-x^\ast|e^{-\lambda^\ast t}\Big(1+\frac{\hat\lambda}{\lambda^\ast}|\bar{x}(0)-x^\ast|^{\hat p-1}\\
&
\times\left(1-e^{-(\hat p-1)\lambda^\ast t}\right)\Big)^{-\frac{1}{\hat p-1}}.
\end{aligned}
$$}
\end{theorem}
\begin{proof}
Consider $V(\bar{x})=\frac{1}{2}|\bar{x}-x^\ast|^2$. From Assumption~1 and the definitions of $\lambda^\ast$ and $\hat\lambda$,
$$
\dot V\leq-\lambda^\ast|\bar{x}-x^\ast|^2-\hat\lambda|\bar{x}-x^\ast|^{\hat p+1}
=-2\lambda^\ast V-2^{\frac{\hat p+1}{2}}\hat\lambda V^{\frac{\hat p+1}{2}}.
$$
For $V>0$, let $W=V^{-\frac{\hat p-1}{2}}$. Then
$$
\dot W\geq(\hat p-1)\lambda^\ast W+(\hat p-1)2^{\frac{\hat p-1}{2}}\hat\lambda.
$$
Consequently,
$$
W(t)\geq e^{(\hat p-1)\lambda^\ast t}\left(W(0)+2^{\frac{\hat p-1}{2}}\frac{\hat\lambda}{\lambda^\ast}\left(1-e^{-(\hat p-1)\lambda^\ast t}\right)\right).
$$
Substituting $W=V^{-\frac{\hat p-1}{2}}$ and $V=\frac{1}{2}|\bar{x}-x^\ast|^2$ gives the stated estimate.
\end{proof}
The estimate in Theorem~1 combines the two decay properties associated with $p^\ast$ and $\hat p$: 
$
\bar{x}(t)-x^\ast|\leq|\bar{x}(0)-x^\ast|e^{-\lambda^\ast t}
$
and
$
|\bar{x}(t)-x^\ast|\leq\left(|\bar{x}(0)-x^\ast|^{1-\hat p}+(\hat p-1)\hat\lambda t\right)^{-\frac{1}{\hat p-1}}.
$
The relation between them can be characterized more precisely as follows.
\begin{corollary}
\textit{Under the conditions of Theorem~1, define
$
\rho^\ast:=(\tfrac{\lambda^\ast}{\hat\lambda})^{\frac{1}{\hat p-1}}.
$
Then the solutions of~\eqref{Lie_sys} satisfy
$$
\begin{aligned}
|\bar{x}(t)-x^\ast|\leq\min\Big\{|&\bar{x}(0)-x^\ast|e^{-\lambda^\ast t},\\
&\left(|\bar{x}(0)-x^\ast|^{1-\hat p}+(\hat p-1)\hat\lambda t\right)^{-\frac{1}{\hat p-1}}\Big\}.
\end{aligned}
$$
Moreover, if $|\bar{x}(0)-x^\ast|\leq\rho^\ast$, then
$$
|\bar{x}(t)-x^\ast|\leq|\bar{x}(0)-x^\ast|e^{-\lambda^\ast t},\qquad t\geq0.
$$
If $|\bar{x}(0)-x^\ast|>\rho^\ast$, there is a
$
\tau^\ast\in\left(0,\frac{(\rho^\ast)^{1-\hat p}-|\bar{x}(0)-x^\ast|^{1-\hat p}}{(\hat p-1)\hat\lambda}\right]
$
such that
$$
|\bar{x}(t)-x^\ast|\leq
\begin{cases}
\begin{aligned}
&\left(|\bar{x}(0)-x^\ast|^{1-\hat p}+(\hat p-1)\hat\lambda t\right)^{-\frac{1}{\hat p-1}}\\
&\hspace{35mm}\text{ for }0\leq t\leq\tau^\ast,
\end{aligned}\\[2mm]
\begin{aligned}
&\rho^\ast e^{-\lambda^\ast(t-\tau^\ast)}\hspace{16mm}\text{ for }t\geq\tau^\ast.
\end{aligned}
\end{cases}
$$}
\end{corollary}
\begin{proof}
The exponential estimate follows directly from Theorem~1 since the second factor in its right-hand side is not greater than one. Furthermore, raising the estimate of Theorem~1 to the power $1-\hat p$ gives
$$
|\bar{x}(t)-x^\ast|^{1-\hat p}\geq|\bar{x}(0)-x^\ast|^{1-\hat p}e^{(\hat p-1)\lambda^\ast t}+\frac{\hat\lambda}{\lambda^\ast}\left(e^{(\hat p-1)\lambda^\ast t}-1\right).
$$
Using $e^s\geq1$ and $e^s-1\geq s$, we obtain
$$
|\bar{x}(t)-x^\ast|\leq\left(|\bar{x}(0)-x^\ast|^{1-\hat p}+(\hat p-1)\hat\lambda t\right)^{-\frac{1}{\hat p-1}}.
$$
The comparison of the two terms follows directly from the definition of $\rho^\ast$. If $|\bar{x}(0)-x^\ast|>\rho^\ast$, asymptotic convergence implies that there exists a first $\tau^\ast>0$ such that $|\bar{x}(\tau^\ast)-x^\ast|=\rho^\ast$. The polynomial estimate gives
$
\tau^\ast\leq\tfrac{(\rho^\ast)^{1-\hat p}-|\bar{x}(0)-x^\ast|^{1-\hat p}}{(\hat p-1)\hat\lambda}.
$
Applying the exponential estimate from $t=\tau^\ast$ gives
$
|\bar{x}(t)-x^\ast|\leq\rho^\ast e^{-\lambda^\ast(t-\tau^\ast)}$ for $ t\geq\tau^\ast.
$
\end{proof}
Thus, solutions starting in $B_{\rho^\ast}(x^\ast)$ admit the exponential decay estimate from the beginning. For solutions starting outside this neighborhood, the polynomial decay estimate associated with $\hat p$ describes the motion until the solution reaches $B_{\rho^\ast}(x^\ast)$, after which the exponential decay estimate applies.
\subsection{Extremum Seeking Design}
We now construct an extremum seeking system whose Lie bracket approximation is given by~\eqref{Lie_sys}.
For each $k\in\{0,\ldots,q\}$, let
$
N_k:=2k+1,
$
and associate with $N_k$ a pair of oscillating inputs $u_{1k}^\omega$, $u_{2k}^\omega$. Consider the extremum seeking system
\begin{equation}
\label{es_sys}
\dot{x}=\sum_{k=0}^{q}\left(g_{1k}(J(x))u_{1k}^\omega(t)+g_{2k}(J(x))u_{2k}^\omega(t)\right).
\end{equation}
We choose
\begin{equation}
\label{es_g}
g_{1k}(z)\equiv1,\qquad g_{2k}(z)=-\gamma_{N_k}z,
\end{equation}
so that
$$
\operatorname{ad}_{g_{1k}}^{N_k}g_{2k}(J(x))=-\gamma_{N_k}J^{(N_k)}(x).
$$
Following~\cite{grushkovskaya2026extremum}, we use the dither signals
\begin{equation}
\label{es_u}
\begin{aligned}
u_{1k}^\omega(t)&=\omega^{\frac{N_k}{N_k+1}}c_{1k}\cos(\kappa_k\omega t),\\
u_{2k}^\omega(t)&=\omega^{\frac{N_k}{N_k+1}}c_{2k}\sin(N_k\kappa_k\omega t),
\end{aligned}
\end{equation}
where $\kappa_k\in\mathbb{N}$ and
$
c_{1k}^{N_k}c_{2k}=(2\kappa_k)^{N_k}N_k!(-1)^k.
$
For each $k$, the corresponding pair of inputs generates the Lie bracket $\operatorname{ad}_{g_{1k}}^{N_k}g_{2k}$. In order to exclude additional Lie brackets resulting from interactions between dithers associated with different orders, we impose a non-resonance condition on the corresponding frequency tuples. 
\begin{assumption}
\textit{The frequency tuples
$
\Omega_k=(\kappa_k,N_k\kappa_k)$,  $k=0,\ldots,q,
$
have no mixed resonances of order up to $N_q+1$. Namely, for any integers $\mu_k,\nu_k$ satisfying
$\displaystyle
0<\sum_{k=0}^{q}(|\mu_k|+|\nu_k|)\leq N_q+1,
$
the relation
$  \displaystyle
\sum_{k=0}^{q}(\mu_k+N_k\nu_k)\kappa_k=0
$
implies
$
\mu_k+N_k\nu_k=0$ for all $ k=0,\ldots,q.
$}
\end{assumption}

\begin{remark}
Assumption~2 is a sufficient non-resonance condition and can be relaxed by taking into account the particular structure of the vector fields. Indeed, $g_{1k}$ are identical for all $k$, while $g_{2k}$ are proportional. Consequently, some mixed resonances correspond to vanishing Lie brackets or to the same derivative directions as those already present in~\eqref{Lie_sys}. We impose Assumption~2 in order to avoid an explicit characterization of such resonances.
\end{remark}

For the terms containing inputs from one pair only, the construction in~\cite{grushkovskaya2026extremum} applies directly: all terms of length not larger than $N_k$ vanish, while the terms of length $N_k+1$ combine into
$
\operatorname{ad}_{g_{1k}}^{N_k}g_{2k}.
$
At the same time, Assumption~2 excludes resonances between inputs corresponding to different pairs.

Since $\kappa_k\in\mathbb{N}$, all the dithers in~\eqref{es_u} are $T$-periodic with
$
T:=\frac{2\pi}{\omega}.
$
Applying the Chen--Fliess representation gives
\begin{equation}
\label{CF_es}
\begin{aligned}
x(T)
&=x^0+T\sum_{k=0}^{q}\operatorname{ad}_{g_{1k}}^{N_k}g_{2k}(J(x^0))+R(T)\\
&=x^0-T\sum_{k=0}^{q}\gamma_{2k+1}J^{(2k+1)}(x^0)+R(T),
\end{aligned}
\end{equation}
where $R(T)$ contains the remaining terms. Moreover, on every compact subset of $D$, there exists  $C_R>0$, independent of $\omega$, such that, for $\omega>1$,
$
|R(T)|\leq C_R\omega^{-1-\frac{1}{N_q+1}}.
$
The above representation indicates that, for sufficiently large $\omega$, the behavior of the solutions of~\eqref{es_sys} over each dither period is close to that of the Lie bracket system~\eqref{Lie_sys}. Therefore, the stability and decay properties established in the previous subsection can be transferred to~\eqref{es_sys} as follows.
\begin{proposition}
\textit{Let Assumptions~1 and~2 hold and assume, in addition, that $J\in C^{2q+2}(D;\mathbb{R})$. Then, for any $\delta>0$ such that $B_\delta(x^\ast)\subset D$, any $\rho\in(0,\delta)$, and any $\lambda\in(0,\lambda^\ast)$, there exists $\bar\omega>0$ such that, for every $\omega\geq\bar\omega$ and $x(0)\in B_\delta(x^\ast)$, the following estimates hold for all $t\geq0$:
\begin{itemize}
\item If $p^\ast=1$, then
$\displaystyle
|x(t)-x^\ast|\leq|x(0)-x^\ast|e^{-\lambda t}+\rho.
$
\item If $p^\ast>1$, then
$$
|x(t)-x^\ast|\leq\left(|x(0)-x^\ast|^{1-p^\ast}+(p^\ast-1)\lambda t\right)^{-\frac{1}{p^\ast-1}}+\rho.
$$
\end{itemize}
Thus, $x^\ast$ is semi-globally practically exponentially stable if $p^\ast=1$, and semi-globally practically asymptotically stable with a polynomial decay rate if $p^\ast>1$.
}
\end{proposition}
\begin{proof}
The proof follows the same arguments as in~\cite{GZE18} and~\cite{grushkovskaya2026extremum}, so we only indicate the main steps which differ due to the representation~\eqref{CF_es}. Fix $\delta>0$ such that $B_\delta(x^\ast)\subset D$, and let $D'\subset D$ be a compact set containing $B_\delta(x^\ast)$ in its interior. For $T=2\pi/\omega$, representation~\eqref{CF_es} holds uniformly on $D'$.

Consider $V(x)=\frac{1}{2}|x-x^\ast|^2$. From Assumption~1 and~\eqref{CF_es},
$$
V(x(T))-V(x^0)\leq
-T\sum_{k=0}^{q}\gamma_{2k+1}\beta_{k1}|x^0-x^\ast|^{p_k+1}
+o(T),
$$
where $o(T)/T\to0$ uniformly on $D'$ as $\omega\to\infty$. Therefore, outside an arbitrary neighborhood of $x^\ast$, for any $\lambda\in(0,\lambda^\ast)$ and sufficiently large $\omega$,
$$
V(x(T))-V(x^0)\leq-T\lambda|x^0-x^\ast|^{p^\ast+1}.
$$
Repeating the period-by-period arguments from~\cite{GZE18} and~\cite{grushkovskaya2026extremum}, and using the estimates of Lemma~1, yields the corresponding practical exponential decay for $p^\ast=1$ and practical polynomial decay for $p^\ast>1$ at the sampling instants $t=nT$.
Furthermore, similarly to~\cite{grushkovskaya2026extremum}, for every $t\in[nT,(n+1)T]$,
$$
|x(t)-x(nT)|\leq C_x\omega^{-\frac{1}{N_q+1}}
$$
with some $C_x>0$ independent of $\omega$. Since $0\leq t-nT<T=2\pi/\omega$, the difference between the decay estimates evaluated at $nT$ and at $t$ also tends to zero uniformly for $x(0)\in B_\delta(x^\ast)$. Hence, by first taking the practical bound at the sampling instants with $\rho/2$ and then choosing $\omega$ sufficiently large, both the inter-sampling deviation and the difference between $nT$ and $t$ can be absorbed into the remaining $\rho/2$. This yields the stated estimates for all $t\geq0$.
\end{proof}
The above proposition characterizes the practical decay rate determined by $p^\ast$. We next show that the combined effect of the terms corresponding to $p^\ast=1$ and $\hat p>1$, established in Theorem~1 for the Lie bracket system, is also preserved by~\eqref{es_sys} up to an    error which can be made arbitrary small by increasing $\omega$.

\begin{theorem}
\label{thm_es_combined}
\textit{Let Assumptions~1 and~2 hold, let $J\in C^{2q+2}(D;\mathbb{R})$, and suppose that $p^\ast=1$ and $\hat p>1$. Then, for any $\delta>0$ such that $B_\delta(x^\ast)\subset D$, any $\rho\in(0,\delta)$, $\lambda\in(0,\lambda^\ast)$, and $\hat\lambda_0\in(0,\hat\lambda)$, there exists $\bar\omega>0$ such that, for every $\omega\geq\bar\omega$ and $x(0)\in B_\delta(x^\ast)$,
$$
\begin{aligned}
|x(t)-x^\ast|\leq&|x(0)-x^\ast|e^{-\lambda t}\Big(1+\frac{\hat\lambda_0}{\lambda}|x(0)-x^\ast|^{\hat p-1}\\
&\times\left(1-e^{-(\hat p-1)\lambda t}\right)\Big)^{-\frac{1}{\hat p-1}}+\rho\text{ for all }t\ge0.
\end{aligned}
$$
}
\end{theorem}

\begin{proof}
The proof follows the same arguments as in~\cite{GZE18} and~\cite{grushkovskaya2026extremum}, so we only indicate the main steps which differ due to the representation~\eqref{CF_es}. Choose $\lambda_1\in(\lambda,\lambda^\ast)$ and $\hat\lambda_1\in(\hat\lambda_0,\hat\lambda)$ such that
$
\frac{\hat\lambda_1}{\lambda_1}>\frac{\hat\lambda_0}{\lambda}.
$
Using~\eqref{CF_es} and Assumption~1, for sufficiently large $\omega$ and outside a neighborhood of $x^\ast$,
$$
V(x(T))-V(x^0)\leq
-T\lambda_1|x^0-x^\ast|^2
-T\hat\lambda_1|x^0-x^\ast|^{\hat p+1},
$$
where $V(x)=\frac{1}{2}|x-x^\ast|^2$.

Let $x_n=x(nT)$ and $V_n=V(x_n)$. Then
$$
V_{n+1}\leq V_n-2\lambda_1TV_n
-2^{\frac{\hat p+1}{2}}\hat\lambda_1T
V_n^{\frac{\hat p+1}{2}}.
$$
For
$
W_n:=V_n^{-\frac{\hat p-1}{2}},
$
and sufficiently large $\omega$, it follows that
$$
W_{n+1}\geq
\left(1+(\hat p-1)\lambda_1T\right)W_n
+(\hat p-1)2^{\frac{\hat p-1}{2}}\hat\lambda_1T.
$$
Consequently,
$$
\begin{aligned}
W_n\geq&
\left(1+(\hat p-1)\lambda_1T\right)^n
\Big(W_0+
2^{\frac{\hat p-1}{2}}\frac{\hat\lambda_1}{\lambda_1}\\
&\times\left(1-
\left(1+(\hat p-1)\lambda_1T\right)^{-n}\right)\Big).
\end{aligned}
$$
Since $T=2\pi/\omega$, for sufficiently large $\omega$,
$$
\left(1+(\hat p-1)\lambda_1T\right)^{\frac{1}{(\hat p-1)T}}
\geq e^\lambda.
$$
Using also $\hat\lambda_1/\lambda_1>\hat\lambda_0/\lambda$, we obtain at the sampling instants
$$
\begin{aligned}
|x(nT)-x^\ast|\leq&|x(0)-x^\ast|e^{-\lambda nT}
\Big(1+\frac{\hat\lambda_0}{\lambda}|x(0)-x^\ast|^{\hat p-1}\\
&\times\left(1-e^{-(\hat p-1)\lambda nT}\right)\Big)^{-\frac{1}{\hat p-1}}
\end{aligned}
$$
until the trajectory reaches an arbitrarily small neighborhood of $x^\ast$.
Furthermore, as in the proof of Proposition~1,
$$
|x(t)-x(nT)|\leq C_x\omega^{-\frac{1}{N_q+1}},
\, t\in[nT,(n+1)T],
$$
and $0\leq t-nT<T$. Since $|x(t)-x(nT)|\to0$ and $t-nT\to0$ uniformly as $\omega\to\infty$, for sufficiently large $\omega$ the estimates obtained at the sampling instants $t=nT$ remain valid for all $t\geq0$ up to the additional term $\rho$.
\end{proof}

Similarly to Corollary~1, Theorem~\ref{thm_es_combined} also yields separate practical estimates corresponding to $p^\ast$ and $\hat p$.
\begin{corollary}
\label{cor_es_combined}
\textit{Under the conditions of Theorem~\ref{thm_es_combined},
$$
\begin{aligned}
|x(t)-x^\ast|\leq&\min\Big\{
|x(0)-x^\ast|e^{-\lambda t},\\
&\left(|x(0)-x^\ast|^{1-\hat p}
+(\hat p-1)\hat\lambda_0t\right)^{-\frac{1}{\hat p-1}}
\Big\}+\rho
\end{aligned}
$$
for all $t\geq0$.
}
\end{corollary}
Thus, the solutions of~\eqref{es_sys} inherit the two decay properties of the Lie bracket system up to the practical error $\rho$: the terms corresponding to $\hat p$ provide a stronger decrease farther from $x^\ast$, while the term corresponding to $p^\ast=1$ gives exponential decay close to $x^\ast$.

An important particular case is given by the family
$$
J_m(x)=(x-x^\ast)^m,\, m\in\{2,4,\ldots,2q+2\}.
$$
For every $2k+1\leq m-1$,
$$
J_m^{(2k+1)}(x)=\frac{m!}{(m-2k-1)!}(x-x^\ast)^{m-2k-1},
$$
while $J_m^{(2k+1)}\equiv0$ for $2k+1>m$.
Hence, independently of the particular value of $m$, the Lie bracket system contains the linear term corresponding to $2k+1=m-1$. In this case,
$
p^\ast=1,\, \hat p=m-1.
$
Thus, the same extremum seeking system~\eqref{es_sys} provides practical exponential decay close to $x^\ast$, while the lower-order derivatives provide higher-power terms farther from the minimum.

\subsection{Extension to the multi-variable case}

The above results can be extended to a class of multivariable cost functions. Let $x\in\mathbb{R}^n$ and, for $k=0,\ldots,q$, define
$$
\mathcal{D}_{2k+1}J(x):=\Big(\frac{\partial^{2k+1}J(x)}{\partial x_1^{2k+1}},\ldots,\frac{\partial^{2k+1}J(x)}{\partial x_n^{2k+1}}\Big)^\top.
$$

\begin{assumption}\label{ass_multi}
The function $J\in C^{2q+1}(D;\mathbb{R})$, $D\subset\mathbb{R}^n$, satisfies the following properties:
\\
A3.1) there exists an $x^\ast\in D$ such that $J(x^\ast)=J^\ast$ and $J(x)>J^\ast$ for all $x\in D\setminus\{x^\ast\}$;
\\
A3.2) for every $k\in\{0,\ldots,q\}$, there exist constants $\beta_{k1},\beta_{k2}\geq0$ and $p_k\geq1$ such that, for all $x\in D$,
$$
(x-x^\ast)^\top\mathcal{D}_{2k+1}J(x)\geq\beta_{k1}\|x-x^\ast\|^{p_k+1},
$$
$$
\|\mathcal{D}_{2k+1}J(x)\|\leq\beta_{k2}\|x-x^\ast\|^{p_k},
$$
and there exists $\bar{k}\in\{0,\ldots,q\}$ such that $\beta_{\bar{k}1}>0$.
\end{assumption}
Consider the Lie bracket system
\begin{equation}
\label{Lie_sys_multi}
\dot{\bar{x}}=-\sum_{k=0}^{q}\gamma_{2k+1}\mathcal{D}_{2k+1}J(\bar{x}),\, \gamma_{2k+1}>0.
\end{equation}
For $k=0$, $\mathcal{D}_1J=\nabla J$, and hence the first term in~\eqref{Lie_sys_multi} coincides with the gradient flow. Define
$$
p^\ast:=\min\{p_k:\beta_{k1}>0\},\, \lambda^\ast:=\sum_{\substack{p_k=p^\ast\\ \beta_{k1}>0}}\gamma_{2k+1}\beta_{k1}.
$$

\begin{proposition}
\label{prop_multi}
Under Assumption~\ref{ass_multi}, the conclusions of Lemma~\ref{lemma_lbs_decay} hold for~\eqref{Lie_sys_multi} with $|\bar{x}-x^\ast|$ replaced by $\|\bar{x}-x^\ast\|$.
\end{proposition}
\begin{proof}
For $V(\bar{x})=\frac{1}{2}\|\bar{x}-x^\ast\|^2$, Assumption~\ref{ass_multi} yields
$$
\begin{aligned}
\dot V&=-\sum_{k=0}^{q}\gamma_{2k+1}(\bar{x}-x^\ast)^\top\mathcal{D}_{2k+1}J(\bar{x})
\leq-\lambda^\ast\|\bar{x}-x^\ast\|^{p^\ast+1}.  
\end{aligned}
$$
The result follows by the same comparison argument as in the proof of Lemma~\ref{lemma_lbs_decay}.
\end{proof}
Note that  Assumption~\ref{ass_multi} is suitable when the relevant derivative fields admit common growth orders. For separable cost functions, however, different coordinates may have different local degrees. This case is covered by the following corollary.
\begin{corollary}
\label{cor_multi_sep}
Let
$
J(x)=\sum_{i=1}^{n}J_i(x_i),
$
where each $J_i\in C^{2q+1}(D_i;\mathbb{R})$ has an isolated minimum at $x_i^\ast$. Assume that, for every $i\in\{1,\ldots,n\}$ and $k\in\{0,\ldots,q\}$, there exist constants $\beta_{ik1},\beta_{ik2}\geq0$ and $p_{ik}\geq1$ such that
$$
\begin{aligned}
&J_i^{(2k+1)}(x_i)(x_i-x_i^\ast)\geq\beta_{ik1}|x_i-x_i^\ast|^{p_{ik}+1},\\
&|J_i^{(2k+1)}(x_i)|\leq\beta_{ik2}|x_i-x_i^\ast|^{p_{ik}},
\end{aligned}
$$
and, for every $i$, there exists at least one $k$ such that $\beta_{ik1}>0$. Then~\eqref{Lie_sys_multi} is decoupled as
$$
\dot{\bar{x}}_i=-\sum_{k=0}^{q}\gamma_{2k+1}J_i^{(2k+1)}(\bar{x}_i),\, i=1,\ldots,n.
$$
Let
$
p_i^\ast:=\min\{p_{ik}:\beta_{ik1}>0\},\, \lambda_i^\ast:=\sum_{\substack{p_{ik}=p_i^\ast\\ \beta_{ik1}>0}}\gamma_{2k+1}\beta_{ik1}.
$
Then, for every $i$, if $p_i^\ast=1$,
$$
|\bar{x}_i(t)-x_i^\ast|\leq|\bar{x}_i(0)-x_i^\ast|e^{-\lambda_i^\ast t},
$$
whereas, if $p_i^\ast>1$,
$$
|\bar{x}_i(t)-x_i^\ast|\leq\Big(|\bar{x}_i(0)-x_i^\ast|^{1-p_i^\ast}+(p_i^\ast-1)\lambda_i^\ast t\Big)^{-\frac{1}{p_i^\ast-1}}.
$$
\end{corollary}
Thus, different coordinates may have different local degrees and corresponding decay rates.

\begin{proof}
For a separable cost, system~\eqref{Lie_sys_multi} decomposes into $n$ scalar systems. The result follows by applying Lemma~\ref{lemma_lbs_decay} to each coordinate.
\end{proof}

The combined decay estimate of Theorem~\ref{thm_lbs_decay} extends in the same way. In particular, defining
$
\hat p:=\max\{p_k:\beta_{k1}>0\},\, \hat\lambda:=\sum_{\substack{p_k=\hat p\\ \beta_{k1}>0}}\gamma_{2k+1}\beta_{k1},
$
the corresponding estimate holds with $|\bar{x}-x^\ast|$ replaced by $\|\bar{x}-x^\ast\|$.

To construct the corresponding extremum seeking system, let $e_i$ denote the $i$th canonical basis vector. A corresponding extremum seeking system can be written as
$$
\dot{x}=\sum_{i=1}^{n}\sum_{k=0}^{q}e_i\Big(u_{1ik}^{\omega}(t)-\gamma_{2k+1}J(x)u_{2ik}^{\omega}(t)\Big).
$$
Assigning to every pair $(i,k)$ dithers of the form~\eqref{es_u}, with frequency tuples satisfying the nonresonance condition jointly over all $(i,k)$, yields~\eqref{Lie_sys_multi}. Indeed, the single-pair terms generate
$
\operatorname{ad}_{e_i}^{2k+1}\Big(-\gamma_{2k+1}J e_i\Big)=-\gamma_{2k+1}\frac{\partial^{2k+1}J}{\partial x_i^{2k+1}}e_i,
$
while additional mixed resonant terms are excluded by the frequency choice. Under the additional regularity $J\in C^{2q+2}(D;\mathbb{R})$, the practical stability results for the extremum seeking system extend analogously.

\section{Numerical Simulations}
\subsection{Single-variable non-polynomial cost function}
We consider the family of non-polynomial cost functions
$$
J_m(x)=2^m\left(
\cosh\Big(\frac{x-x^\ast}{2}\right)
-\sum_{\ell=0}^{m/2-1}
\frac{(x-x^\ast)^{2\ell}}{2^{2\ell}(2\ell)!}
\Big),
$$
with $ m\in\{2,4,6\}$ and $x^\ast=0$. By construction,
$$
J_m(x)=\frac{(x-x^\ast)^m}{m!}
+O\left(|x-x^\ast|^{m+2}\right)\text{ as } x\to x^\ast,
$$
so that $m$ determines the local degree of the cost, while $J_m$ is not a polynomial. For all three cost functions, we use the same gains
$
\gamma_1=\gamma_3=\gamma_5=1,
$
and  compare three extremum seeking systems of form~\eqref{es_sys} whose corresponding Lie bracket systems are
\begin{align}
\dot{\bar{x}}&=-J^{(1)}(\bar{x}),
\label{lbs_num_1}\\
\dot{\bar{x}}&=-\left(J^{(1)}(\bar{x})+J^{(3)}(\bar{x})\right),
\label{lbs_num_13}\\
\dot{\bar{x}}&=-\left(J^{(1)}(\bar{x})+J^{(3)}(\bar{x})+J^{(5)}(\bar{x})\right),
\label{lbs_num_135}
\end{align}
respectively.
The corresponding frequency parameters are chosen as $\kappa_0=1$, $(\kappa_0,\kappa_1)=(1,4)$, and $(\kappa_0,\kappa_1,\kappa_2)=(1,13,22)$, respectively. In all simulations, $\omega=10^4$ and $x(0)=2$.

It is important to note that, in contrast to the polynomial case, the higher-order derivatives do not vanish even for the lowest considered degree $m=2$. In this case, all three Lie bracket systems have a linear leading term near $x^\ast$, while the additional higher-order brackets increase its coefficient:
$$
\begin{aligned}
&J_2^{(1)}(x)=(x-x^\ast)+O(|x-x^\ast|^3),\\
&J_2^{(3)}(x)=\frac{1}{4}(x-x^\ast)+O(|x-x^\ast|^3),\\
&J_2^{(5)}(x)=\frac{1}{16}(x-x^\ast)+O(|x-x^\ast|^3).
\end{aligned}
$$
For $m=4$, the local behavior of the relevant derivatives is
$$
J_4^{(1)}(x)=\frac{(x-x^\ast)^3}{3!}
+O\left(|x-x^\ast|^5\right),
$$
$$
\begin{aligned}
&J_4^{(3)}(x)=(x-x^\ast)+O(|x-x^\ast|^3),\\
&J_4^{(5)}(x)=\frac{1}{4}(x-x^\ast)+O(|x-x^\ast|^3).
\end{aligned}
$$
Hence, the first-order design has polynomial local decay, whereas the addition of the third-order Lie bracket introduces a linear term and yields exponential local decay. The fifth-order bracket remains nonzero and further contributes to the linear part.

Finally, for $m=6$,
$$
\begin{aligned}
&J_6^{(1)}(x)=\frac{(x-x^\ast)^5}{5!}+O(|x-x^\ast|^7),\\
&J_6^{(3)}(x)=\frac{(x-x^\ast)^3}{3!}+O(|x-x^\ast|^5),\\
&J_6^{(5)}(x)=(x-x^\ast)+O(|x-x^\ast|^3).
\end{aligned}
$$
Thus, the three designs  introduce fifth-, third-, and first-power terms in the local Lie bracket dynamics. In particular, the fifth-order bracket provides the linear term responsible for exponential decay close to the minimizer.

\begin{figure*}
\begin{center}
        \includegraphics[width=0.77\linewidth]{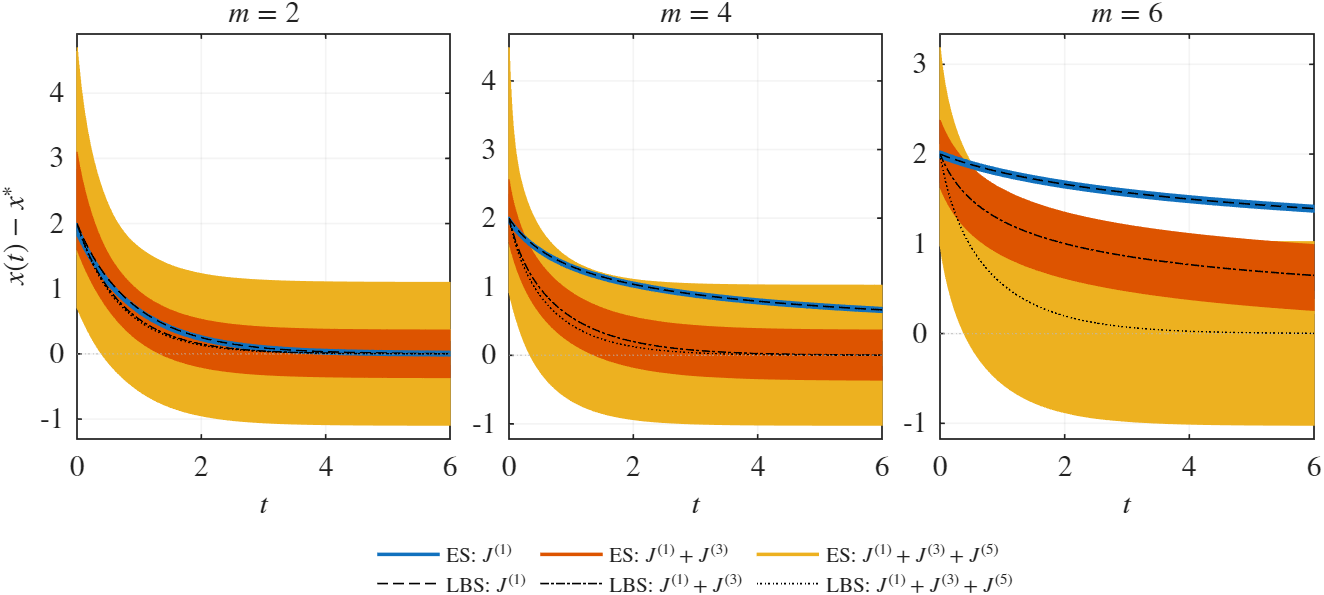}
\end{center}
\caption{Trajectories of the three extremum seeking systems for the non-polynomial costs $J_m$. The blue, orange, and yellow ES trajectories correspond to the Lie bracket systems~\eqref{lbs_num_1}, \eqref{lbs_num_13}, and~\eqref{lbs_num_135}, respectively. The left, middle, and right plots correspond to $m=2$, $m=4$, and $m=6$.}\label{fig_ex1}
\end{figure*}

The trajectories are shown in Fig.~\ref{fig_ex1}. The example illustrates that the proposed design is not restricted to polynomial cost functions. The three costs have different local degrees, while the same controller gains are used in all cases and no information on the local degree is used in their selection. Moreover, for a non-polynomial cost, derivatives of order higher than the local degree need not vanish and may provide additional descent directions. The simulations also show the improvement in the local decay obtained by including higher-order Lie brackets, together with the finite-frequency oscillations of the corresponding extremum seeking systems.

\subsection{Multi-variable non-homogeneous cost function}
To illustrate the multivariable case with different local degrees, consider
$$
J(x)=\frac{(x_1^2+x_2^2)^2}{4!}+\frac{x_3^2}{2!},
$$
with $x^\ast=(0,0,0)^\top$. The variables $x_1$ and $x_2$ form a coupled block of local degree four, whereas the local degree with respect to $x_3$ is two. The relevant derivative vectors are
$$
\mathcal{D}_1J(x)=\Big(\frac{x_1(x_1^2+x_2^2)}{3!},\,\frac{x_2(x_1^2+x_2^2)}{3!},\,x_3\Big)^\top,
$$
and
$\displaystyle
\mathcal{D}_3J(x)=(x_1,\,x_2,\,0)^\top.
$
We compare the extremum seeking systems corresponding to the Lie bracket systems
$$
\dot{\bar{x}}=-\mathcal{D}_1J(\bar{x})
$$
and
$$
\dot{\bar{x}}=-\mathcal{D}_1J(\bar{x})-\mathcal{D}_3J(\bar{x}),
$$
using $\gamma_1=\gamma_3=1$. In the first design, the dynamics of $\bar{x}_1$ and $\bar{x}_2$ contain only cubic terms, while the $\bar{x}_3$-dynamics is linear. Adding the third-order Lie bracket introduces the linear terms $-\bar{x}_1$ and $-\bar{x}_2$, while leaving the $\bar{x}_3$-dynamics unchanged since $\partial^3J/\partial x_3^3=0$.

\begin{figure*}[t]
    \centering
    \includegraphics[width=0.8\textwidth]{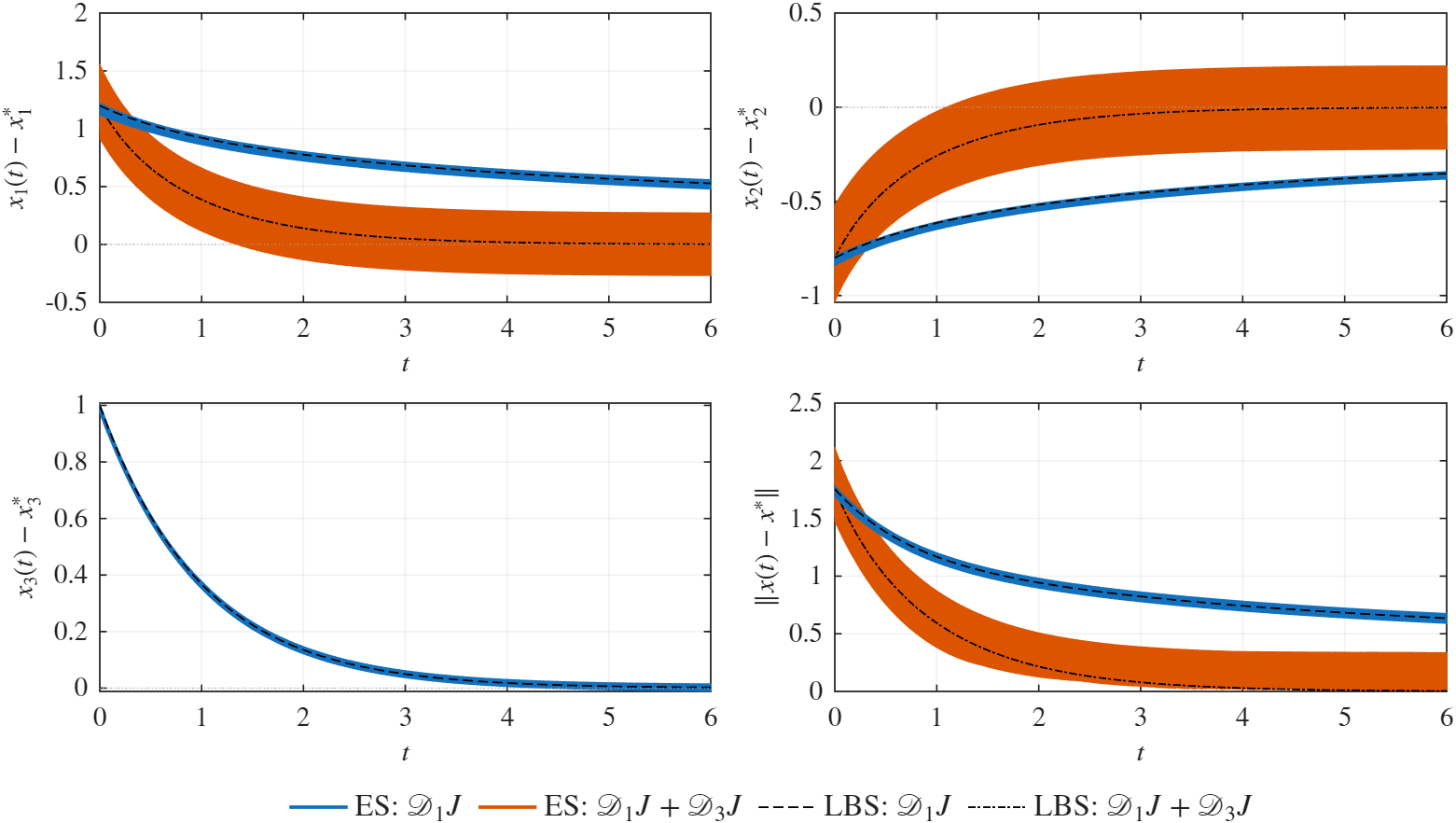}
    \caption{Trajectories of the multivariable extremum seeking systems for $J(x)=(x_1^2+x_2^2)^2/4!+x_3^2/2!$ and the corresponding Lie bracket systems. The blue curves correspond to $\mathcal{D}_1J$, while the orange curves correspond to $\mathcal{D}_1J+\mathcal{D}_3J$.}
    \label{fig:exmpl2}
\end{figure*}

The trajectories are shown in Fig.~\ref{fig:exmpl2}. As expected, the two designs exhibit essentially the same behavior in the $x_3$-coordinate, whereas the addition of the third-order Lie bracket substantially improves the decay of $x_1$ and $x_2$. Consequently, the norm $\|x(t)-x^\ast\|$ also decreases considerably faster for the combined design. This example illustrates that different directions may have different local degrees and that the additional higher-order Lie brackets become active only in those directions where they provide a nonzero descent term.
\section{Conclusion and Future Work}
In this paper, we introduce new results that built on, and explore further, the recent theory and framework on constructing higher-order Lie bracket approximations of control-affine systems and extremum seeking (ES) \cite{pokhrel2026higher,grushkovskaya2026extremum}. We provide results that show the possibility of exciting first-order Lie-bracket flow combined with higher-order Lie-bracket flow, allowing for an accelerated new ES design. That is, the proposed design is based on Lie bracket dynamics that combine the gradient term with higher-order derivative terms, benefiting from the stronger polynomial-type decrease away from the minimizer and exponential decay near the minimizer whenever a linear-like higher-order term is present. The provided numerical simulations, including
multi-variable and non-polynomial flat-bottomed objective
functions, demonstrate the potential of the proposed combined design.

It is important to highlight that the proposed design supports and sheds the light on future possibility where accelerated and exponentially convergent ES becomes possible even if the local degree of the objective function is unknown (i.e., $m$ is totally unknown) since we demonstrated successfully and rigorously that Lie brackets of different orders can be excited simultaneously. It is also of crucial importance to improve the proposed ES by addressing the large oscillation problem present in this design, similar to what was noted in \cite{pokhrel2026higher,grushkovskaya2026extremum}.
\bibliographystyle{ieeetr}
\bibliography{biblio_ES,references}
\end{document}